\documentclass[12pt,a4paper]{amsart}

\usepackage {amsfonts}
\usepackage[english]{babel}
\def\b{\beta}
\def\g{\gamma}
\def\G{\Gamma}
\def\d{\delta}
\def\a{\alpha}
\def\p{\varphi}
\def\e{\varepsilon}
\def\l{\lambda}

\def\cS{{\mathcal S}}
\def\la{\langle}
\def\ra{\rangle}
\def\R{{\mathbb R}}
\def\C{{\mathbb C}}
\def\Z{{\mathbb Z}}

\def\bs{~\hfill\rule{7pt}{7pt}}

\DeclareMathOperator{\Lin}{Lin}
\DeclareMathOperator{\supp}{supp}
\DeclareMathOperator{\dist}{dist}

\newtheorem{Th}{Theorem}
\newtheorem*{Cor}{Corollary}
\newtheorem{Lem}{Lemma}

\begin{document}

\title{Local Versions of the Wiener--L\'{e}vy Theorem}

\author{S.Yu.Favorov}

\address{Serhii Favorov,
\newline\hphantom{iii}  Karazin's Kharkiv National University
\newline\hphantom{iii} Svobody sq., 4,
\newline\hphantom{iii} 61022, Kharkiv, Ukraine}
\email{sfavorov@gmail.com}

\maketitle {\small
\begin{quote}
\noindent{\bf Abstract.}
Let h be a real-analytic function in the neighborhood of some compact set K on the plane. We show that for any complex measure on the Euclidean space of a finite total variation without singular components with the Fourier--Stieltjes transform f(y) there exists another measure of a finite total variation with the  Fourier  transform g(y) with the property  g(y)=h(f(y)) for each y such that f(y) belongs to K.
\medskip


\end{quote}
}

\medskip
{\bf 1. Introduction}. It is well known  that for each absolutely convergent Fourier series $F(t)$ such that $F(t)\neq0$ for all $t$ the function $1/F(t)$ also has an absolutely convergent Fourier-series expansion (the Wiener Theorem).
Its natural generalization is known as the Wiener--L\'evy Theorem (see, for example, \cite{Z}, Ch.VI):
\begin{Th}
Let
$$
F(t)=\sum_{n\in\Z}c_ne^{2\pi int}
$$
 be an absolutely convergent Fourier series,  and $h(z)$  be a holomorphic function on a neighborhood of the closure of the range of  $F$. Then the function $h(F(t))$ admits  an absolutely convergent Fourier series expansion as well.
 \end{Th}
 Clearly, for $h(z)=1/z$ we get the Wiener Theorem.
 \medskip

 The next variant of this theorem for functions on $\R$ is also known as the Wiener--L\'evy Theorem (see for example \cite{A}, Ch.III or \cite{R}, Ch.I)
 \begin{Th}\label{2}
Let $\hat f(y)$ be the Fourier transform of some function $f\in L^1(\R)$, and $h(z)$  be a holomorphic function on a neighborhood of the closure of the range of $\hat f$ such that $h(0)=0$. Then there is a function
$g\in L^1(\R)$ such that its Fourier transform $\hat g$ coincides with $h(\hat f(y))$.
 \end{Th}

This theorem admits a generalization to functions from $L^1(G)$, where $G$ is an arbitrary locally compact abelian group, and one can replace a holomorphic function $h$ by any real-analytic function. Then if we replace here the absolutely continuous measure $f(x)dx$ by a pure point measure $\sum_n a_n\d_{\l_n}$ ($\d_s$, as usual, means the unit mass at the point $s$) with $\sum_n |a_n|<\infty$, we obtain the Wiener--L\'evy Theorem for Diriclet series. On the other hand, for each locally compact abelian non discrete group $G$ there exists a measure $\mu$ with a finite total variation such that values of its Fourier  transform $\hat\mu$ is bounded away from zero on the dual group  $\hat G$ (hence the function $1/z$ is holomorphic on the closure of the set $\hat\mu(\hat G)$), but there is no measure $\nu$ on $G$ such that its Fourier  transform $\hat\nu=1/\hat\mu$. Also, the requirement of the real-analyticity is  necessary for most groups $G$ (in particular for $G=\R^d$) to fulfill the Wiener--L\'evy Theorem (see \cite{R1}, Ch.5, 6).
\medskip

Also note that the Fourier transforms of pure point measures  were considered earlier in \cite{Fe} in connection with the problem of complete reconstruction of  band-limited functions.
\medskip

 The local form of the Wiener--L\'evy Theorem is of greatest interest for our study (see \cite{R}, Ch.6):
 \begin{Th}\label{3}
Let $G$ be a locally compact abelian group, let $K$ be a compact subset of the dual group $\hat G$, let $f\in L^1(G)$, and let $h(z)$  be a holomorphic function on a neighborhood of the closure of the set $\hat f(K)$. Then there is a function
$g\in L^1(G)$ such that its Fourier transform $\hat g$ coincides with $h(\hat f(y))$ for all $y\in K$.
 \end{Th}
In our article we consider the case when $K$ is a compact subset of the complex plane $\C$ and the function $h(z)$ is analytic (or real-analytic) on a neighborhood of $K$. Of course if $\hat G$ is a compact abelian group then $\hat f^{-1}(K)$ is a compact subset of $\hat G$, hence we are in the conditions of Theorem \ref{3}. But for $G=\R^d$ we obtain a result  stronger than Theorem \ref{3}. Theorems of this type were used by us to study Poisson measures in \cite{F} -- \cite{F2}.
\medskip

{\bf 2. Notations and Preliminaries}. To formulate our results we have to recall some definitions.

Denote by $M(G)$ the set of complex measures on the locally compact group $G$ with a finite total variation $\|\mu\|$, by $M_d(G)$ the set of pure point measures from $M(G)$, and by $M_{ad}$ the set of measures from $M(G)$ containing only pure point and absolutely continuous (with respect to the Haar measure) components, i.e., without  singular components. The Fourier  transform of $\mu\in M(G)$ is defined by the equality
 $$
    \hat\mu(y)=\int_G (-x,y)\mu(dx),\qquad y\in\hat G,
 $$
where $\hat G$ is the group of characters on $G$, and $(x,y)$ means the action of the character $y$ on $x\in G$. In particular, in the case $G=\R^d$ we have
$$
    \hat\mu(y)=\int_{\R^d}e^{-2\pi i\la x,y\ra}\mu(dx),\qquad y\in\R^d,
 $$
where $\la x,y\ra$ means the scalar product of $x$ and $y$. If the measure $\mu$ is absolute continuous and $\mu(dx)=f(x)dx$, we will also write
$$
    \hat f(y)=\int_{\R^d}f(x)e^{-2\pi i\la x,y\ra}dx.
 $$
Furthermore, a complex-valued function $h$, defined on an open set $V\subset\C$ is said to be real-analytic on $V$ if to every point $z\in V$ there corresponds the expansion
\begin{equation}\label{ra}
   h(\xi+i\eta)=\sum_{k,n=0}^\infty c_{k,n}(\xi-\Re z)^k(\eta-\Im z)^n,\quad \xi,\,\eta\in\R, \quad c_{k,n}\in\C,
\end{equation}
which converges in some disc
$$
D(z)=\{(\xi,\eta)\in\R^2:\, |\xi-\Re z|^2+|\eta-\Im z|^2<r_z^2\}.
$$
Note that series (\ref{ra}) converges also in the ball
$$
B(z)=\{(\xi,\eta)\in\C^2:\, |\xi-\Re z|^2+|\eta-\Im z|^2<r_z^2\},
$$
and for any intersecting balls $B(z_1)$ and $B(z_2)$ the set
$$
B(z_1)\cap B(z_2)\cap\R^2=D(z_1)\cap D(z_2)
$$
 is the set of uniqueness for analytic functions of two variables. Therefore, if the function $h(z)$ is real-analytic in a neighborhood of some compact set $K\subset\C$, then it has a continuation to the neighborhood
 $\cup_{z\in K}B(z)\subset\C^2$ of $K$ as an analytic function of two complex variables $\xi,\,\eta$.
\newpage

{\bf 3. Main results}. \begin{Th}\label{4}
Let $\mu$ be a measure from $M_{ad}(\R^d)$, let $h(z)$ be a real-analytic function on a neighborhood of some compact set $K\subset\C$. Then there is a measure $\nu\in M_{ad}(\R^d)$ such that, for every $y\in\R^d$ for which $\hat\mu(y)\in K$, we have $\hat\nu(y)=h(\hat\mu(y))$.
\end{Th}
In particular, if  $h(z)=1/z$ or $h(z)=1/|z|^\a$ for $|z|\ge\e$ and $h(z)=0$ for $|z|\le\e/2$, we obtain the following result:
\begin{Cor}
For any $\mu\in M_{ad}(\R^d)$  and $\e>0,\ \a>0$ there are measures $\nu_\e,\nu_{\a,\e}\in M_{ad}(\R^d)$  such that in the case $|\hat\mu(y)|\ge\e$ we have $\hat\nu_\e(y)=1/\hat\mu(y)$, $\hat\nu_{\a,\e}(y)=1/|\hat\mu(y)|^\a$,   and in the case $|\hat\mu(y)|\le\e/2$ we have $\hat\nu_\e(y)=\hat\nu_{\a,\e}(y)=0$.
\end{Cor}
The reasoning in the proof of Theorem \ref{4} also provides  the following statement:
\begin{Th}\label{5}
Let $\mu$ be a measure from $M_d(G)$ for a locally compact abelian group $G$, and let $h(z)$ be a real-analytic function on a neighborhood of some compact set $K\subset\C$. Then  there is a measure $\nu\in M_d(G)$  such that, for every $y\in\hat G$ for which $\hat\mu(y)\in K$, we have $\hat\nu(y)=h(\hat\mu(y))$. The support of the measure $\nu$ lies in $\Lin_\Z\supp\mu$.
\end{Th}
Here $\supp\mu$  means the set $\{x\in G: \mu(\{x\})\neq0\}$ if  $\mu\in M_d(G)$.
\medskip

For the case $G=\R^d$ and holomorphic $h$ on a neighborhood $V\subset\C$ of $K$ Theorem \ref{5} was proved in \cite{F}. Its application allowed us to obtain in \cite{F} some strengthening of one of the theorems in the theory of Fourier quasicrystals. Another application of Theorem \ref{5} to Kahane's property of discrete sets see \cite{F1}.
\medskip

{\bf 4. Auxiliary lemmas and their proofs}.
We will use   Schwartz' space $\cS(\R^d)$ of rapidly decreasing $C^\infty$-functions on $\R^d$ with the topology defined by a countable number of norms
 $$
   N_n(\p)=\sup_{x\in\R^d}\left\{(1+|x|)^n\max_{k_1+\dots+k_d\le n}\left|\frac{\partial^{k_1}}{\partial x_1^{k_1}}\dots\frac{\partial^{k_d}}{\partial x_d^{k_d}}\p(x)\right|\right\},\ n=0,1,2,\dots
 $$
The Fourier transform is a continuous linear one-to-one mapping of $\cS(\R^d)$ onto $\cS(\R^d)$, and the set of $C^\infty$-functions with bounded support is dense in $\cS(\R^d)$ (see \cite{R2}).
 \begin{Lem}\label{L1}
For every $f\in L^1(\R^d)$ and every $\e>0$ there is $v\in\cS(\R^d)$ such that $\|f-v\|_{L^1}<\e$ and $\hat v$ has a compact support.
\end{Lem}
{\sl Proof}. Take $f_1\in L^1(\R^d)$ such that $\|f-f_1\|_{L^1}<\e/3$ and $f_1$ has a compact support. The convolution $f_2=f_1\star\p$ with a suitable $C^\infty$-function $\p(x)\ge0$ with support in a small ball such that $\int_{\R^d}\p(x)dx=1$ has the properties $\|f_2-f_1\|_{L^1}<\e/3$ and $f_2\in\cS(\R^d)$. Therefore, $\hat f_2\in\cS(R^d)$ as well, and there is a sequence  of $C^\infty$-functions with compact supports that converges to $\hat f_2$ in the space $\cS(R^d)$. Let $\{v_n\}$ be the images the functions from this sequence under the inverse Fourier transform. Clearly, $v_n\to f_2$ in the space $\cS(R^d)$, therefore,
 $$
   \|f_2-v_n\|_{L^1}\le  \sup_{\R^d}(1+|x|)^{d+1}|f_2(x)-v_n(x)|\int_{\R^d}(1+|x|)^{-d-1}dx\le C(d)\, N_{d+1}(f_2-v_n)\to 0
  $$
 as $n\to\infty$. Hence, $\|f_2-v_n\|_{L^1}<\e/3$ for a suitable $v_n$, and $\|f-v_n\|_{L^1}<\e$.  \bs

\begin{Lem}\label{L2}
There is a constant $C=C(r,d)$ such that for every $v\in\cS(\R^d)$ with the property $\supp\hat v\subset B(0,r)$ we get
$$
\|v\|_{L^1}\le C(r,d)\left(\|\hat v\|_\infty+\sum_{j=1}^d\|(\partial^l/\partial y_j^l)\hat v\|_\infty\right),
$$
with $l=d+1$ for odd $d$ and $l=d+2$ for even $d$.
\end{Lem}
{\sl Proof}. We have
 $$
  \hat v(y)=\int_{\R^d}v(x)e^{-2\pi i\la x,y\ra} dx,\quad (\partial^l/\partial y_j^l)\hat v(y)=(-2\pi i)^l\int_{\R^d}x_j^l v(x)e^{-2\pi i\la x,y\ra} dx,\ j=1,\dots,d.
 $$
 Hence, $v(x)(1+\sum_{j=1}^d |x_j|^l)$ is the inverse Fourier transform of the function
 $$
 \hat v(y)+\left(\frac{i}{2\pi}\right)^l\sum_{j=1}^d (\partial^l/\partial y_j^l)\hat v(y).
 $$
Since $\supp\hat v\subset B(0,r)$, we get
$$
\|v\|_{L^1}\le C\left(\|\hat v\|_\infty+\sum_{j=1}^d\|(\partial^l/\partial y_j^l)\hat v\|_\infty\right) \int_{\R^d}\frac{dx}{1+|x_1|^l+|x_2|^l+\dots+|x_d|^l}
$$
with a constant $C$ depending on $d$ and $r$. \bs
\begin{Lem}\label{L3}
Let $T(\Theta,\tau)$ be $C^\infty$-function in variables $\Theta=(\theta_1,\dots,\theta_N)\in \R^N$ and $\tau\in [0,1]^2$, and let $T(\Theta,\tau)$ be periodic with periods $1$ in each coordinate  $\theta_1,\dots,\theta_N$.
Then its Fourier series
$$
   T(\Theta,\tau)=\sum_{k\in\Z^N}b_k(\tau)e^{2\pi i\la k,\Theta\ra},\qquad k=(k_1,\dots,k_N),
$$
is absolutely convergent and $\sum_k|b_k(\tau)|<C$ uniformly in $\tau$.
\end{Lem}
{\sl Proof}. We have
 $$
   b_k(\tau)=\int_{[0,1]^N}T(\Theta,\tau)e^{-2\pi i\la k,\Theta\ra}d\Theta.
 $$
   Integrating this equality twice in parts over each variable $\theta_j$ such that $j\in J(k)=\{j:\,k_j\neq0\}$, we get
 $$
   |b_k(\tau)|\le \sup_{(\Theta,\tau)\in[0,1]^{N+2}}\left|\left(\prod_{j\in J(k)}(-4\pi^2 k_j^2)^{-1}\partial^2/\partial\theta_j^2\right)T(\Theta,\tau)\right|.
 $$
 Taking into account that every derivative of $T(\Theta,\tau)$ is uniformly bounded in $\tau\in [0,1]^2$, we get the estimate
$$
   |b_k(\tau)|\le C\min\{1,k_1^{-2}\}\cdots\min\{1,k_N^{-2}\},
 $$
where $C$ depends on neither $\tau$ nor $k$. This estimate implies the assertion of the lemma. \bs
\medskip

{\bf 5. Proofs of the main results}.
\medskip

{\sl Proof of Theorem \ref{4}}. Let $U$ be an open set in $\C^2$ such that $h$ has a holomorphic continuation to $U$ as a function of two complex variables. Set $\e<(1/13)\dist(K,\partial U)$.
  Let $\p(|z|)$ be $C^\infty$-differentiable nonnegative function with a support in $B(0,\e)\subset\C^2$ such that $\int_{\C^2}\p(|\zeta|)m(d\zeta)=1$ (here $m(d\zeta)$ means the Lebesgue measure in $\C^2$). Consider
  $C^\infty$-function
$$
     H(z)=\int_{\dist(\zeta,K)<9\e} h(\zeta)\p(|z-\zeta|)m(d\zeta).
$$
 If $\dist(z,K)<8\e$, we get
 $$
     H(z)=\int_{|\zeta|\le\e} h(z-\zeta)\p(|\zeta|)m(d\zeta).
$$
  Since an average over any sphere of a holomorphic function of many variables equals the meaning of the function in the center of the sphere, we obtain that $H(z)=h(z)$
  on the set $\{z:\dist(z,K)<7\e\}$ and $H(z)=0$ on the set $\{z:\dist(z,K)>10\e\}$.

     Let
  $$
  \mu=fdx+\sum_n a_n\d_{\g_n},\quad f\in L^1(\R^d),\quad \sum_n|a_n|<\infty.
  $$
  Using Lemma \ref{L1}, take a function $v\in\cS(\R^d)$ such that $\|f-v\|_{L^1}<\e$ and $\supp\hat v$ is a compact set. Pick $N<\infty$ such that $\sum_{n>N}|a_n|<\e$, and define the measure $s=\sum_{n=1}^N a_n\d_{\g_n}$.
  Note that
  $$
  \hat s(y)=\sum_{n\le N}a_n e^{-2\pi i\la\g_n,y\ra}.
  $$
   Since $\|\mu-vdx-s\|<2\e$, we see that $\|\hat\mu(y)-\hat v(y)-\hat s(y)\|_\infty<2\e$.
  Put
  $$
  \a(y)=\Re(\hat v(y)+\hat s(y)),\ \b(y)=\Im(\hat v(y)+\hat s(y)).
  $$
  Consider the function
  $$
     F(y)=\frac{1}{(2\pi i)^2}\int_{|\a(y)-\zeta_1|=3\e}\int_{|\b(y)-\zeta_2|=3\e}\frac{H(\zeta_1+i\zeta_2)d\zeta_1 d\zeta_2}{(\zeta_1-\Re\hat\mu(y))(\zeta_2-\Im\hat\mu(y))}.
  $$
If $\hat\mu(y)\in K$, then $\dist(\a(y)+i\b(y),K)<2\e$. Therefore,
$$
    \mathcal E=\{(\zeta_1,\zeta_2):\,|\zeta_1-\a(y)|\le 3\e,\,|\zeta_2-\b(y)|\le 3\e\}\subset\{z:\dist(z,K)<7\e\},
$$
and $H(z)=h(z)$ in a neighborhood of $\mathcal E$. Using the Cauchy integral formula for the polydisk $\mathcal E$, we obtain
  $$
       F(y)=h(\Re\hat\mu(y)+i\Im\hat\mu(y))=h(\hat\mu(y)).
  $$

Furthermore, we have for all $y\in\R^d$
  $$
     F(y)=\int_0^1\int_0^1\frac{H(\a(y)+3\e e^{3\pi i\tau_1}+i\b(y)+i3\e e^{2\pi i\tau_2})9\e^2 e^{2\pi i(\tau_1+\tau_2)}}{(\a(y)+3\e e^{3\pi i\tau_1}-\Re\hat\mu(y))(\b(y)+3\e e^{3\pi i\tau_2}-\Im\hat\mu(y))}d\tau_1\,d\tau_2.
  $$
 Since
 $$
 \left|\frac{\Re\hat\mu(y)-\a(y)}{3\e e^{2\pi i\tau_1}}\right|<2/3,\qquad \left|\frac{\Im\hat\mu(y)-\b(y)}{3\e e^{2\pi i\tau_2}}\right|<2/3,
 $$
  we get
 $$
  \frac{1}{\left(1-\frac{\Re\hat\mu(y)-\a(y)}{3\e e^{2\pi i\tau_1}}\right)\left(1-\frac{\Im\hat\mu(y)-\b(y)}{3\e e^{2\pi i\tau_2}}\right)}= \sum_{p,q=0}^\infty\left(\frac{\Re\hat\mu(y)-\a(y)}{3\e e^{2\pi i\tau_1}}\right)^p
  \left(\frac{\Im\hat\mu(y)-\b(y)}{3\e e^{2\pi i\tau_2}}\right)^q,
 $$
 and
$$
     F(y)= \sum_{p,q=0}^\infty\left(\frac{\Re\hat\mu(y)-\a(y)}{3\e}\right)^p
  \left(\frac{\Im\hat\mu(y)-\b(y)}{3\e}\right)^q
  \int_0^1\int_0^1 \frac{A(y,\tau_1,\tau_2)+D(y,\tau_1,\tau_2)}{e^{2\pi i(p\tau_1+q\tau_2)}}d\tau_1\,d\tau_2,
  $$
with
 $$
  A(y,\tau_1,\tau_2)=H(\hat v(y)+\hat s(y)+3\e e^{2\pi i\tau_1}+i3\e e^{2\pi i\tau_2})-H(\hat s(y)+3\e e^{2\pi i\tau_1}+i3\e e^{2\pi i\tau_2}),
 $$
$$
  D(y,\tau_1,\tau_2)=H(\hat s(y)+3\e e^{2\pi i\tau_1}+i3\e e^{2\pi i\tau_2}).
 $$
 Define two measures on $\R^d$
 $$
 \l_R(x)=1/2\left(\mu(x)-v(x)dx-s(x)+\overline{\mu(-x)-v(-x)dx-s(-x)}\right),
 $$
 $$
\l_I(x)=(1/2i)\left(\mu(x)-v(x)dx-s(x)-\overline{\mu(-x)-v(-x)dx-s(-x)}\right).
$$
 It is easily seen that $\|\l_R\|<2\e,\  \|\l_I\|<2\e$, and
 $$
 \hat\l_R(y)=\Re\hat\mu(y)-\a(y),\ \hat\l_I(y)=\Im\hat\mu(y)-\b(y).
$$
 Since the Fourier transform of convolution of measures equals the product of the Fourier transform of the measures, we get
  \begin{equation}\label{con}
   \left[(\l_R/3\e)^{*p}*(\l_I/3\e)^{*q}\right]^{\widehat{\quad}}= \left(\frac{\Re\hat\mu(y)-\a(y)}{3\e}\right)^p
  \left(\frac{\Im\hat\mu(y)-\b(y)}{3\e}\right)^q.
  \end{equation}
Also, the variation of convolution of measures does not exceed the product of variations of the measures, hence
 \begin{equation}\label{conv}
 \|(\l_R/3\e)^{*p}*(\l_I/3\e)^{*q}\|<(2/3)^{p+q}.
 \end{equation}

 On the other hand, since  $\supp A(y,\tau_1,\tau_2)\subset\supp\hat v$ and $A(y,\tau_1,\tau_2)$ is $C^\infty$ function, we see that $A(y,\tau_1,\tau_2)\in\cS(\R^d)$. Therefore there exists
 $u_{\tau_1,\tau_2}(x)\in\cS(\R^d)$ such that $\hat u_{\tau_1,\tau_2}(y)=A(y,\tau_1,\tau_2)$ for every fixed $\tau_1,\tau_2$. Then the function $A(y,\tau_1,\tau_2)$ and all its derivatives of order
 at most $d+2$ are bounded uniformly in $\tau_1,\,\tau_2\in [0,1]^2$. By Lemma \ref{L2}, $\|u_{\tau_1,\tau_2}\|_{L^1}$ is uniformly bounded too. Set
  $$
     \kappa_{p,q}(x)=\int_0^1\int_0^1 \frac{u_{\tau_1,\tau_2}(x)}{e^{2\pi i(p\tau_1+q\tau_2)}}d\tau_1\,d\tau_2\in L^1(\R^d).
  $$
  By Fubini's Theorem,
  \begin{equation}\label{norm}
  \sup_{p,q}\|\kappa_{p,q}\|_{L^1}<\infty,
  \end{equation}
  and
  \begin{equation}\label{nor}
     \hat\kappa_{p,q}(y)=\int_0^1\int_0^1 \frac{A(y,\tau_1,\tau_2)}{e^{2\pi i(p\tau_1+q\tau_2)}}d\tau_1\,d\tau_2.
  \end{equation}

    Next, apply Lemma \ref{L3} to the function $H\left(\sum_{n\le N}a_n e^{2\pi i\theta_n}+3\e e^{2\pi i\tau_1}+i3\e e^{2\pi i\tau_2}\right)$. We get
 \begin{equation}\label{pd}
  H\left(\sum_{n\le N}a_n e^{2\pi i\theta_n}+3\e e^{2\pi i\tau_1}+i3\e e^{2\pi i\tau_2}\right)=\sum_{k\in\Z^N}b_k(\tau_1,\tau_2)e^{2\pi i\la k,\Theta\ra},
\end{equation}
with the condition
 \begin{equation}\label{s}
  \sup_{\tau_1,\tau_2}\sum_k|b_k(\tau_1,\tau_2)|<\infty.
 \end{equation}
 If we replace in (\ref{pd}) $\theta_n$ by $-\la\g_n,y\ra$  for each $n$, we get the function\footnote{Note that some $\rho_k$ may coincide}
 $$
 \sum_{k\in\Z}b_k(\tau_1,\tau_2)e^{2\pi i\la\rho_k,y\ra},\qquad \rho_k\in\Lin_\Z\{\g_n\}_{n=1}^\infty.
 $$
  This function is the Fourier transform of the measure
 $\sum_{k\in\Z}b_k(\tau_1,\tau_2)\d_{-\rho_k}$. Set
 $$
 \nu_{p,q}=\sum_{k\in\Z}c_k(p,q)\d_{-\rho_k}\quad\mbox{with}\quad c_k(p,q)=\int_0^1\int_0^1 \frac{b_k(\tau_1,\tau_2)}{e^{2\pi i(p\tau_1+q\tau_2)}}d\tau_1\,d\tau_2.
 $$
 It follows from (\ref{pd}), (\ref{s}) and  Fubini's Theorem that
 \begin{equation}\label{meas}
   \sup_{p,q}\|\nu_{p,q}\|<\infty,
 \end{equation}
 and
\begin{equation}\label{mea}
  \hat\nu_{p,q}(y)= \int_0^1\int_0^1 \frac{D(y,\tau_1,\tau_2)}{e^{2\pi i(p\tau_1+q\tau_2)}}d\tau_1\,d\tau_2.
\end{equation}

 Finally put
 \begin{equation}\label{sum}
      \nu= \sum_{p,q=0}^\infty(\l_R/3\e)^{*p}*(\l_I/3\e)^{*q}*(\kappa_{p,q}dx+\nu_{p,q}).
 \end{equation}
 We have
  $$
     \|\nu\|\le \sum_{p,q=0}^\infty\|\l_R/3\e\|^p\|\l_I/3\e\|^q(\|\kappa_{p,q}\|_{L^1}+\|\nu_{p,q}\|).
  $$
 It follows from (\ref{conv}), (\ref{norm}), and (\ref{meas}), that $\nu$
 has a finite total variation, and, by (\ref{con}), (\ref{nor}), and (\ref{mea}), that $\hat\nu(y)=F(y)$. \bs

\medskip
{\sl Proof of Theorem \ref{5}}. Let $\mu=\sum_n a_n\d_{\g_n}$ with $\g_n\in G$ and $\sum_n|a_n|<\infty$. Then
$$
\hat\mu(y)=\sum_n a_n(-\g_n,y), \qquad y\in\hat G.
 $$
 Replace $e^{-2\pi i\la y,\g_n\ra}$ by $(-\g_n,y)$ in the previous proof. Further, we do not use Lemmas \ref{L1} and \ref{L2}, but put  $A(y,\tau_1,\tau_2)\equiv0, \ u(x,\tau_1,\tau_2)\equiv0$, and $k_{p,q}(x)\equiv0\ \forall p,q$. Then, repeating the reasoning in the proof of Theorem \ref{4}, we obtain the assertion of Theorem \ref{5}. \bs

 \medskip
{\bf Remark}. If the function $h$ is holomorphic in a neighborhood of the compact set $K$, then there is an alternative proof of Theorem \ref{5}. Indeed, let $\G$ denote the group $G$ with respect to the discrete topology. Clearly, every pure point measure $\mu\in M(G)$ is the function $f\in L^1(\G)$  at the same time. Therefore $\hat\mu$ extends to the continuous function $\hat f$ on the compact group $\hat\G$. Then $\hat f^{-1}(K)$ is a compact subset of $\hat\G$, and we may apply Theorem \ref{3}. In order to obtain a statement about the support of the measure $\nu$, one has to replace the group $\G$ by the group $\Lin_\Z\supp\mu$.

\bigskip

I am very grateful to Professor Hans Georg Feichtinger for pointing me the papers, which contain results close to those in my article, and for his attention to my work.

\end{document}